\documentclass{article}

\usepackage{PRIMEarxiv}

\usepackage[utf8]{inputenc} 
\usepackage[T1]{fontenc}    
\usepackage{hyperref}       
\usepackage{url}            
\usepackage{booktabs}       
\usepackage{amsfonts}       
\usepackage{nicefrac}       
\usepackage{microtype}      
\usepackage{lipsum}
\usepackage{fancyhdr}       
\usepackage{graphicx}       
\graphicspath{{media/}}     
\usepackage{amsmath}
\usepackage{amsthm}
\usepackage{fdsymbol}
\usepackage{float}
\usepackage{tikz}

\theoremstyle{plain} 
\newtheorem{theorem}{Theorem}[section]
\theoremstyle{definition}
\newtheorem{corollary}[theorem]{Corollary}
\newtheorem{lemma}[theorem]{Lemma}
\newtheorem{proposition}{Proposition}[section]
\newtheorem{definition}{Definition}[section]

\title{Constructions of $A_\alpha$-cospectral graphs using some corona operations
}

\author{
  Najiya V K, Chithra A V \\
 Department of Mathematics \\
 National Institute of Technology, Calicut\\
 Kerala, India-673601\\
  \texttt{najiya\_p190046ma@nitc.ac.in, chithra@nitc.ac.in} \\
}

\begin{document}
\maketitle

\begin{abstract}
Let $ G_1 \circledast G_2$,$ G_1 \sqcupdot G_2 $ and $ G_1 \sqcupplus G_2$ denote the total corona, $Q$-vertex corona and $Q$-edge corona of two graphs $ G_1$ and $ G_2 $, respectively. In this paper, we compute the $A_\alpha$-spectrum of $ G_1 \circledast G_2$,$ G_1 \sqcupdot G_2 $ and $ G_1 \sqcupplus G_2$ for regular graphs $ G_1$ and $ G_2$. As an application, we construct infinitely many pairs of $A_\alpha$-cospectral graphs.
\end{abstract}


\section{Introduction}
Throughout the manuscript, we focus only on simple undirected graphs.  Let $ G$ denote a graph with vertex set $V( G)=\{u_{1}, u_{2},...,u_{n}\}$ and edge set $E( G)=\{e_{1}, e_{2},...,e_{m}\}$. The adjacency matrix $A( G)$ of $ G$ is a square symmetric matrix of order $n$ whose $(i,j)^{th}$ entry is 1 if the vertices $u_{i}$ and $u_{j}$  are adjacent and is 0 otherwise. The degree matrix $D( G)$ is a diagonal matrix with diagonal entries corresponding to vertex degrees. The $A_\alpha$-matrix\cite{nikiforov2017merging} $A_\alpha ( G)$ of $ G$ is defined as
$A_\alpha ( G) = \alpha D( G) + (1 - \alpha)A( G)$ for any real $\alpha \in [0, 1]$. Specifically, for $\alpha=0, \frac{1}{2}$ and $1$, $A_\alpha( G)$ becomes $A( G)$, $\frac{1}{2}Q( G)$ and $D( G)$, respectively. 

The characteristic polynomial of the matrix $M$ of order $n$ is defined as $\phi(M,\lambda)=|\lambda I_n-M|$, where $I_n$ is the identity matrix of order $n$. The roots for this polynomial are called the $M$-eigenvalues.
Let $\lambda_1(A_\alpha ( G)) \geq \lambda_2(A_\alpha ( G)) \geq \dots \geq
\lambda_n (A_\alpha ( G))$ represent the $A_\alpha$-eigenvalues of $ G$. The set of all $A_\alpha ( G)$ eigenvalues along with their multiplicities is called the $A_\alpha$-spectrum of $ G$. For an r-regular graph $ G$, $\lambda_i(A_\alpha( G))=\alpha r+(1-\alpha)\lambda_i(A( G))$, $i=1,2,...,n$. Two graphs are said to be $A_\alpha$-cospectral if they have the same $A_\alpha$-spectrum. Otherwise, they are non $A_\alpha$-cospectral graphs. $A_\alpha$- cospectral graphs are always $A$-cospectral graphs, but $A$-cospectral graphs need not always be $A_\alpha$-cospectral. However, regular $A$-cospectral graphs are $A_\alpha$-cospectral. 

The adjacency energy\cite{bapat2010graphs} $\varepsilon( G)$ of $ G$ is defined as the sum of the absolute values of eigenvalues of $ G$. The $A_\alpha$-energy\cite{pirzada2021alpha} $\varepsilon_\alpha( G)$ of $ G$ is defined as $\displaystyle\sum_{i=1}^n\left|\lambda_i(A_\alpha( G))-\frac{2\alpha m}{n}\right|$. For a regular graph $ G$, $\varepsilon_\alpha( G)=(1-\alpha)\varepsilon( G)$.  
 
 The incidence matrix $R( G)$ of a graph $ G$ is the matrix of order $n\times m$, whose $(i,j)^{th}$ entry is 1 if $u_i$ is incident to $e_j$ and 0 otherwise. From \cite{cvetkovic1980spectra},  $R( G)R( G)^T=A( G)+D( G)$ and $R( G)^TR( G)=B( G)+2I_m$, where $B$ is the adjacency matrix of line graph of $ G$. For an $r$-regular graph $ G$, $R( G)R( G)^T=A( G)+rI_n$ . Throughout this paper, $I_n$ is the identity matrix of order $n$, $O_{m\times n}$ is the $0$ matrix of order $m\times n$, $J_n$ is the matrix of order $n$ with all entries equal to one and $J_{s\times t}$ is the matrix of order $s\times t$ with all entries equal to one. Let $A$ and $B$ be matrices of order $a \times b$ and $c \times d$, respectively, then the Kronecker product of $A$ and $B$, denoted $A\otimes B$, is the $ac\times bd$ block matrix $[a_{ij} B]$. 

 The corona\cite{frucht1970corona} of two graphs $ G_1$, $ G_2$, denoted by $ G_1\circ  G_2$, is constructed by taking one copy of $ G_1$ and $n_1$ copies of $ G_2$ and joining the ith vertex of $ G_1$ to every vertex in the ith copy of $ G_2$ for each $i$. The $Q$ graph of $ G$, $Q( G)$, is the graph obtained from $ G$ by inserting a new vertex into every edge of $ G$ and then joining pairs of these new vertices that lie on adjacent edges of $ G$. The total graph $T( G)$ of a graph $ G$ is the graph whose set of vertices is the union of the set of vertices and of the set of edges of $ G$, with two vertices of $T( G)$ being adjacent if and only if the corresponding elements of $ G$ are adjacent or incident. Let $ G$ be a graph with vertex set $V ( G) = \{v_1 , v_2 , \dots , v_n \}$. Consider a set $U ( G) = \{u_1 , u_2 ,\dots , u_n \}$ of $n$ vertices. Make $u_i$ adjacent to all the vertices in the neighbuorhood of $v_i $, for each $i = 1, 2, \dots , n$. The graph thus obtained is called the splitting graph of $ G$, denoted by $Splt( G)$. Let $ G_1$ and $ G_2$ be two graphs with $n_1$ and $n_2$ vertices, respectively. The neighbourhood corona, $ G_1 \star  G_2 $, of $ G_1$ and $ G_2$ is obtained by taking $n_1$ copies of $ G_2$ and joining all the vertices of the ith copy of $ G_2$ with the vertices adjacent to $v_i$ in $ G_1 $, for $i = 1, 2,\dots , n_1$.
 In this paper, we find the $A_\alpha$-spectrum of a few corona products of graphs. 
  
The rest of the paper is organised as follows. In Section \ref{prelims}, we present some definitions and lemmas used to prove our results. In Sections \ref{totcor}, \ref{splcor}, and \ref{qcor}, we compute the $A_\alpha$-spectrum of the total corona, $Q$-vertex corona, and $Q$-edge corona of two graphs. Also, we construct an infinite number of pairs of $A_\alpha$-cospectral graphs that are non-regular. 

\section{Preliminaries}\label{prelims}

 \begin{definition}\cite{zhu2016spectra}
     The total corona of $ G_1$ and $ G_2$ , denoted by $  G_1 \circledast  G_2 $, is obtained by taking one copy of $T (  G_1 )$ and $n_1$ copies of $ G_2$, and joining the $i$th vertex of $ G_1$ to every vertex in the $i$th copy of $ G_2 $.
 \end{definition}

 \begin{definition}\cite{thomas2023spectrum}
     The splitting corona, $ G_1 \boxplus  G_2$ of $ G_1$ and $ G_2$ is obtained from $Spl( G_1 )$ and $n_1$ copies of $ G_2$ by joining every vertex of the ith copy of $ G_2$ to the ith vertex of $V ( G_1 )$, for $i = 1, 2,\dots, n_1 $.
 \end{definition}

 \begin{definition}\cite{thomas2023spectrum}
     The splitting add vertex corona, $ G_1 \boxtimes  G_2$ of $ G_1$ and $ G_2$ is obtained from $Spl ( G_1 )$ and $n_1$ copies of $ G_2$ by joining every vertex of the ith copy of $ G_2$ to the ith vertex of $U ( G_1 )$, for $i = 1, 2, \dots , n_1 $.
 \end{definition}

 \begin{definition}\cite{thomas2023spectrum}
     The splitting neighbourhood corona, $ G_1 \divideontimes  G_2$ of $ G_1$ and $ G_2$ is obtained from $Splt ( G_1 )$ and $n_1$ copies of $ G_2$ by joining every vertex of the ith copy of $ G_2$ to the vertices adjacent to the ith vertex in $V ( G_1 )$, for $i = 1, 2, \dots , n_1 $.
 \end{definition}


Now, we define new graph operations based on $Q$-graphs as follows.

\begin{definition}
    The $Q$-vertex corona of $ G_1$ and $ G_2$, denoted by $ G_1 \sqcupdot  G_2$, is the graph obtained from $Q\left( G_1\right)$ and $n_1$ copies of $ G_2$ by joining the $i$th vertex of $V\left( G_1\right)$ to every vertex in the $i$th copy of $ G_2$.
\end{definition}

\begin{definition}
    The $Q$-edge corona of $ G_1$ and $ G_2$, denoted by $ G_1 \sqcupplus  G_2$, is the graph obtained from $Q\left( G_1\right)$ and $m_1$ copies of $ G_2$ by joining the $i$th vertex of $E( G_1)$ to every vertex in the $i$th copy of $ G_2$.
\end{definition}

\begin{figure}[H]
\begin{center}
\begin{tikzpicture}[scale=1,auto=center] 
 \tikzset{dark/.style={circle,fill=black}}
 \tikzset{red/.style={circle,fill=red}}
 \tikzset{white/.style={circle,draw=white}}
  
  \node [red] (a1) at (0,0)  {};  
  \node [red] (a2) at (1,0)  {};  
  \node [red] (a3) at (1,1)  {};  
  \node [red] (a4) at (0,1){};
  
  \node [dark] (a5) at (0.5,-1) {};
  \node [dark] (a6) at (2,0.5) {};
   \node [dark] (a7) at (0.5,2){} ;  
  \node [dark] (a8) at (-1,0.5)  {};

  \node [dark] (a9) at (0,-2)  {};  
  \node [dark] (a10) at (1,-2)  {}; 
  \node [dark] (a11) at (3,0)  {};  
  \node [dark] (a12) at (3,1)  {};  
  \node [dark] (a13) at (1,3)  {};  
  \node [dark] (a14) at (0,3){};
   \node [dark] (a15) at (-2,0)  {};  
  \node [dark] (a16) at (-2,1){};

  \draw (a1) -- (a2);
  \draw (a2) -- (a3);  
  \draw (a3) -- (a4);  
  \draw (a4) -- (a1);  
  
  \draw (a5) -- (a1);
  \draw (a5) -- (a2);  
  \draw (a6) -- (a2);  
  \draw (a6) -- (a3); 
  \draw (a7) -- (a3);
  \draw (a7) -- (a4);
  \draw (a8) -- (a4);
  \draw (a8) -- (a1);

   \draw (a10) -- (a9);  
  \draw (a11) -- (a12);  
  \draw (a13) -- (a14);
  \draw (a15) -- (a16); 
  
  \draw (a5) -- (a9);  
  \draw (a5) -- (a10);  
  \draw (a6) -- (a11);
  \draw (a6) -- (a12);  
  \draw (a7) -- (a13);  
  \draw (a7) -- (a14);  
  \draw (a8) -- (a15);  
  \draw (a8) -- (a16);  
  
\end{tikzpicture}  
\end{center}
\caption{$C_4\sqcupdot K_2$} \label{qvecor}
\end{figure}

\begin{figure}[H]
\begin{center}
\begin{tikzpicture}[scale=1,auto=center] 
 \tikzset{dark/.style={circle,fill=black}}
 \tikzset{red/.style={circle,fill=red}}
 \tikzset{white/.style={circle,draw=white}}
  
  \node [dark] (a1) at (0,0)  {};  
  \node [dark] (a2) at (1,0)  {};  
  \node [dark] (a3) at (1,1)  {};  
  \node [dark] (a4) at (0,1){};
  
  \node [red] (a5) at (0.5,-1) {};
  \node [red] (a6) at (2,0.5) {};
   \node [red] (a7) at (0.5,2){} ;  
  \node [red] (a8) at (-1,0.5)  {};

  \node [dark] (a9) at (0,-2)  {};  
  \node [dark] (a10) at (1,-2)  {}; 
  \node [dark] (a11) at (3,0)  {};  
  \node [dark] (a12) at (3,1)  {};  
  \node [dark] (a13) at (1,3)  {};  
  \node [dark] (a14) at (0,3){};
   \node [dark] (a15) at (-2,0)  {};  
  \node [dark] (a16) at (-2,1){};

  \draw (a5) -- (a1);
  \draw (a5) -- (a2);  
  \draw (a6) -- (a2);  
  \draw (a6) -- (a3);  
  \draw (a7) -- (a3);
  \draw (a7) -- (a4);  
  \draw (a8) -- (a4);  
  \draw (a8) -- (a1); 

  \draw (a5) -- (a6);
  \draw (a6) -- (a7);
  \draw (a7) -- (a8);
  \draw (a8) -- (a5);

   \draw (a10) -- (a9);  
  \draw (a11) -- (a12);  
  \draw (a13) -- (a14);
  \draw (a15) -- (a16); 
  
  \draw (a5) -- (a9);  
  \draw (a5) -- (a10);  
  \draw (a6) -- (a11);
  \draw (a6) -- (a12);  
  \draw (a7) -- (a13);  
  \draw (a7) -- (a14);  
  \draw (a8) -- (a15);  
  \draw (a8) -- (a16);  
  
\end{tikzpicture}  
\end{center}
\caption{$C_4\sqcupplus K_2$} \label{qedcor}
\end{figure}
 
 \begin{lemma}\cite{cvetkovic1980spectra}\label{uvwx}
	Let $P,Q,R$ and $S$ be matrices and $$M=\begin{pmatrix}
	P&Q\\
	R&S
	\end{pmatrix}. 
	$$
	
	If $P$ is invertible, then $\det(M)=\det(P)\det(S-RP^{-1}Q)$.
	\medskip
	
	If $S$ is invertible, then $\det(M)=\det(S)\det(P-QS^{-1}R).$ 
	\medskip
	
	If $P$ and $R$ commute, then $\det(M)=\det(PS-QR)$.
\end{lemma}

\begin{definition}\cite{cui2012spectrum}\label{gammad}
	The $M$-coronal $\Gamma_M(\lambda)$  of a matrix $M$ of order $n$ is defined as the sum of the entries of the matrix $(\lambda I_n-M)^{-1}$ (if it exists), that is, $$\Gamma_M(\lambda)=J_{n\times 1}^T(\lambda I_n-M)^{-1}J_{n\times 1}.$$ \\
	If each row sum of an $n\times n$ matrix $M$ is constant, say a, then 
	$\Gamma_{M}(\lambda)=\frac{n}{\lambda-a}$.
\end{definition}
\begin{lemma}\cite{mcleman2011spectra}\label{gamma}
	Let $ G$ be an $r$-regular graph on $n$  vertices,  then $\Gamma_{A( G)}(\lambda)=\frac{n}{\lambda-r}$.
\end{lemma}

\begin{theorem}\cite{barik2018spectra}
    If $ G$ is an $r$-regular graph on $n$ vertices, then the adjacency spectrum of line graph of $ G$ is $$\left(\begin{array}{cc}
        -2 & \lambda_i+r-2 \\
        m-n & 1
    \end{array}\right),$$ where $\lambda_i$'s are the adjacency eigenvalues of $ G$.
\end{theorem}

\section{\texorpdfstring{$A_\alpha$}{A alpha}-spectrum of total corona}\label{totcor}
In this section, we compute of $A_\alpha$-spectrum of $ G_1\circledast  G_2$, the total corona of the graphs $ G_1$ and $ G_2 $. First, we obtain an expression for the $A_\alpha$-characteristic polynomial of $ G_1 \circledast G_2$, where $ G_1$ is an $r_1$-reugular graph and $ G_2$ is an arbitrary graph.

\begin{proposition}
     Let $ G_1$ be an $r_1$-regular graph on $n_1$ vertices and $ G_2$ be an arbitrary graph on $n_2$ vertices. Then
     \begin{align*}
         \phi(A_\alpha  \left( G_1 \circledast  G_2\right),\lambda) =&\left(\lambda-2(\alpha r_1-1+\alpha)\right)^{m_1-n_1}\prod_{i=1}^{n_2}\left(\lambda-\alpha-\lambda_i(A_{\alpha_2})\right)^{n_1} \\
&\hspace{-2cm}\prod_{i=1}^{n_1}\left((\lambda-2(\alpha r_1-1+\alpha))(\lambda-(3r_1+n_2)\alpha-(1-\alpha)^2\Gamma_{A_{\alpha_2}}(\lambda-\alpha)+r_1)\right.\\
&\hspace{-2cm}\left.-(1-\alpha)(\lambda_i(A_1)+r_1)(2\lambda-(5r_1+n_2)\alpha+3(1-\alpha)+r_1-(1-\alpha)^2\Gamma_{A_{\alpha_2}}(\lambda-\alpha))+(1-\alpha)^2(\lambda_i(A_1)+r_1)^2\right).
     \end{align*}
\end{proposition}

\begin{proof}
    With a suitable ordering of the vertices of $ G_1 \circledast  G_2$ , we get

\begin{align*}
A_\alpha\left( G_1 \circledast  G_2\right)=&\left(\begin{array}{ccc}
(r_1+n_2)\alpha I_{n_1}+A_{\alpha_1} & (1-\alpha)R_1 & (1-\alpha)I_{n_1}\otimes J_{1\times n_2} \\
(1-\alpha)R_1^T & 2\alpha r_1I_{m_1}+(1-\alpha)B_1 & O \\
(1-\alpha)I_{n_1}\otimes J_{n_2 \times 1}  & O & I_{n_1}\otimes (\alpha I_{n_2}+A_{\alpha_2})
\end{array}\right).\\
\left|\lambda I-A_\alpha\left( G_1 \circledast  G_2\right)\right|=&\left|\begin{array}{ccc}
(\lambda-(r_1+n_2)\alpha) I_{n_1}-A_{\alpha_1} & -(1-\alpha)R_1 & -(1-\alpha)I_{n_1}\otimes J_{1\times n_2} \\
-(1-\alpha)R_1^T & (\lambda-2\alpha r_1)I_{m_1}-(1-\alpha)B_1 & O \\
-(1-\alpha)I_{n_1}\otimes J_{n_2 \times 1}  & O & I_{n_1}\otimes ((\lambda-\alpha) I_{n_2}-A_{\alpha_2})
\end{array}\right|\\
=&\prod_{i=1}^{n_2}(\lambda-\alpha-\lambda_i(A_{\alpha_2}))^{n_1}\det S\\
\intertext{where }
\det S=&\left|\begin{array}{cc}
    (\lambda-(r_1+n_2)\alpha-(1-\alpha)^2\Gamma_{A_{\alpha_2}}(\lambda-\alpha)) I_{n_1}-A_{\alpha_1} & -(1-\alpha)R_1 \\
    -(1-\alpha)R_1^T & (\lambda-2\alpha r_1)I_{m_1}-(1-\alpha)B_1
\end{array}\right|.
\intertext{Now performing row transformations on rows $\mathcal{R}_1$ and $\mathcal{R}_2$ by $\mathcal{R}_2\rightarrow R_1^T\mathcal{R}_1+\mathcal{R}_2$ followed by $\mathcal{R}_1\rightarrow\mathcal{R}_1+\frac{(1-\alpha)R_1\mathcal{R}_2}{\lambda-2\alpha r_1+2(1-\alpha)}$ we get}
&\hspace{-3.5cm}\det S=\left(\lambda-2(\alpha r_1-1+\alpha)\right)^{m_1-n_1}\prod_{i=1}^{n_1}\left((\lambda-2(\alpha r_1-1+\alpha))(\lambda-(3r_1+n_2)\alpha-(1-\alpha)^2\Gamma_{A_{\alpha_2}}(\lambda-\alpha)+r_1)\right.\\
&\hspace{-2.5cm}\left.-(1-\alpha)(\lambda_i(A_1)+r_1)(2\lambda-(5r_1+n_2)\alpha+3(1-\alpha)+r_1-(1-\alpha)^2\Gamma_{A_{\alpha_2}}(\lambda-\alpha))+(1-\alpha)^2(\lambda_i(A_1)+r_1)^2\right).
\end{align*}
Then the result follows.
\end{proof}

Now, in the following corollary, we obtain the $A_\alpha$-eigenvalues of $ G_1 \circledast  G_2$, where both $ G_1$ and $ G_2$ are regular graphs.

\begin{corollary}
    Let $ G_i$ be $r_i$-regular graph with $n_i$ vertices for $i=1,2$. Then the $A_\alpha$-spectrum of $ G_1 \circledast  G_2$ consists of:
    \begin{enumerate}
    \item $2(\alpha r_1-1+\alpha)$ repeated $m_1-n_1$ times,
        \item $\alpha+\lambda_i(A_{\alpha_2})$ repeated $n_1$ times for each $i=2,3,\dots,n_2$,
        \item three roots of the equation
        $(\lambda-2(\alpha r_1-1+\alpha))((\lambda-\alpha-r_2)(\lambda-(3r_1+n_2)\alpha+r_1)-(1-\alpha)^2n_2)-(1-\alpha)(\lambda_i(A_1)+r_1)((\lambda-\alpha-r_2)(2\lambda-(5r_1+n_2)\alpha+3(1-\alpha)+r_1)-(1-\alpha)^2n_2)+(1-\alpha)^2(\lambda-\alpha-r_2)(\lambda_i(A_1)+r_1)^2=0,$ for $i=1,2,\dots,n_1$.
    \end{enumerate}
\end{corollary}

Finally, as an application, we construct new pairs of $A_\alpha$-cospectral graphs in the following corollary.

\begin{corollary}
    \begin{enumerate}
    \item Let $ G_1$ and $ G_2$ be two $A$-cospectral regular graphs and $H$ be an arbitrary graph. Then the graphs $ G_1\circledast H$ and $ G_2\circledast H$ are $A_\alpha$-cospectral.
    \item Let $H_1$ and $H_2$ be two $A_\alpha$-cospectral graphs with $\Gamma_{A_\alpha(H_1)}(\lambda) = \Gamma_{A_\alpha(H_2)}(\lambda)$ for $\alpha \in [0, 1]$. If $ G$ is a regular graph, then the graphs $ G\circledast H_1$ and $ G\circledast H_2$ are $A_\alpha$-cospectral.
\end{enumerate}
\end{corollary}

\section{\texorpdfstring{$A_\alpha$}{A alpha}-spectrum of corona based on splitting graph}\label{splcor}
In this section, we compute of $A_\alpha$-spectrum of $ G_1\boxplus  G_2$,$ G_1\boxtimes  G_2$ and $ G_1\divideontimes  G_2$ of the graphs $ G_1$ and $ G_2 $. First, we obtain an expression for the $A_\alpha$-characteristic polynomial of these operations, where $ G_1$ is an $r_1$-reugular graph and $ G_2$ is an arbitrary graph.

\begin{proposition}\label{splitcor}
     Let $ G_1$ be an $r_1$-regular graph on $n_1$ vertices and $ G_2$ be an arbitrary graph on $n_2$ vertices. Then
     \begin{align*}
         \phi(A_\alpha  \left( G_1 \boxplus  G_2\right),\lambda) =&\prod_{i=1}^{n_2}\left(\lambda-\alpha-\lambda_i(A_{\alpha_2})\right)^{n_1} \\
&\hspace{-2cm}\prod_{i=1}^{n_1}(\lambda-\alpha r_1)(\lambda-(2r_1+n_2)\alpha-(1-\alpha)^2\Gamma_{A_{\alpha_2}}(\lambda-\alpha)-(1-\alpha)\lambda_i(A_1))-(1-\alpha)^2\lambda_i(A_1).
     \end{align*}
\end{proposition}

\begin{proof}
    With a suitable ordering of the vertices of $ G_1 \boxplus  G_2$ , we get

\begin{align*}
A_\alpha\left( G_1 \boxplus  G_2\right)=&\left(\begin{array}{ccc}
(2r_1+n_2)\alpha I_{n_1}+(1-\alpha)A_1 & (1-\alpha)A_1 & (1-\alpha)I_{n_1}\otimes J_{1\times n_2} \\
(1-\alpha)A_1 & \alpha r_1I_{n_1} & O \\
(1-\alpha)I_{n_1}\otimes J_{n_2 \times 1}  & O & I_{n_1}\otimes (\alpha I_{n_2}+A_{\alpha_2})
\end{array}\right).\\
\left|\lambda I-A_\alpha\left( G_1 \circledast  G_2\right)\right|=&\left|\begin{array}{ccc}
(\lambda-(2r_1+n_2)\alpha) I_{n_1}-(1-\alpha)A_1 & -(1-\alpha)A_1 & -(1-\alpha)I_{n_1}\otimes J_{1\times n_2} \\
-(1-\alpha)A_1 & (\lambda-\alpha r_1)I_{n_1} & O \\
-(1-\alpha)I_{n_1}\otimes J_{n_2 \times 1}  & O & I_{n_1}\otimes ((\lambda-\alpha) I_{n_2}-A_{\alpha_2})
\end{array}\right|\\
=&\prod_{i=1}^{n_2}(\lambda-\alpha-\lambda_i(A_{\alpha_2}))^{n_1}\det S\\
\intertext{where }
\det S=&\left|\begin{array}{cc}
    (\lambda-(2r_1+n_2)\alpha) I_{n_1}-(1-\alpha)A_1-(1-\alpha)^2\Gamma_{A_{\alpha_2}}(\lambda-\alpha)I_{n_1} & -(1-\alpha)A_1  \\
-(1-\alpha)A_1 & (\lambda-\alpha r_1)I_{n_1} 
\end{array}\right|.\\
\intertext{Then by Lemma \ref{uvwx}}
=&\prod_{i=1}^{n_1}(\lambda-\alpha r_1)(\lambda-(2r_1+n_2)\alpha-(1-\alpha)^2\Gamma_{A_{\alpha_2}}(\lambda-\alpha)-(1-\alpha)\lambda_i(A_1))-(1-\alpha)^2\lambda_i(A_1).
\end{align*}
Then the result follows.
\end{proof}

Now, in the following corollary, we obtain the $A_\alpha$-eigenvalues of $ G_1 \boxplus  G_2$, where both $ G_1$ and $ G_2$ are regular graphs.

\begin{corollary}
    Let $ G_i$ be $r_i$-regular graph with $n_i$ vertices for $i=1,2$. Then the $A_\alpha$-spectrum of $ G_1 \boxplus  G_2$ consists of:
    \begin{enumerate}
        \item $\alpha+\lambda_i(A_{\alpha_2})$ repeated $n_1$ times for each $i=2,3,\dots,n_2$,
        \item two roots of the equation
       $(\lambda-\alpha r_1)((\lambda-\alpha-r_2)(\lambda-(2r_1+n_2)\alpha-(1-\alpha)\lambda_i(A_1))-(1-\alpha)^2n_2)-(1-\alpha)^2(\lambda-\alpha-r_2)\lambda_i(A_1)=0,$ for $i=1,2,\dots,n_1$.
    \end{enumerate}
\end{corollary}

Finally, as an application, we construct new pairs of $A_\alpha$-cospectral graphs in the following corollary.

\begin{corollary}
    \begin{enumerate}
    \item Let $ G_1$ and $ G_2$ be two $A$-cospectral regular graphs and $H$ be an arbitrary graph. Then the graphs $ G_1\boxplus H$ and $ G_2\boxplus H$ are $A_\alpha$-cospectral.
    \item Let $H_1$ and $H_2$ be two $A_\alpha$-cospectral graphs with $\Gamma_{A_\alpha(H_1)}(\lambda) = \Gamma_{A_\alpha(H_2)}(\lambda)$ for $\alpha \in [0, 1]$. If $ G$ is a regular graph, then the graphs $ G\boxplus H_1$ and $ G\boxplus H_2$ are $A_\alpha$-cospectral.
\end{enumerate}
\end{corollary}

Next, we compute $A_\alpha$-spectrum of splitting add vertex corona of $ G_1$ and $ G_2$, where $ G_1$ is an $r_1$-reugular graph and $ G_2$ is an arbitrary graph.

\begin{proposition}
     Let $ G_1$ be an $r_1$-regular graph on $n_1$ vertices and $ G_2$ be an arbitrary graph on $n_2$ vertices. Then
     \begin{align*}
         \phi(A_\alpha  \left( G_1 \boxtimes  G_2\right),\lambda) =&\prod_{i=1}^{n_2}\left(\lambda-\alpha-\lambda_i(A_{\alpha_2})\right)^{n_1} \\
&\hspace{-2cm}\prod_{i=1}^{n_1}(\lambda-2\alpha r_1-(1-\alpha) \lambda_i(A_1))(\lambda-(r_1+n_2)\alpha-(1-\alpha)^2\Gamma_{A_{\alpha_2}}(\lambda-\alpha))-(1-\alpha)^2\lambda_i(A_1)^2.
     \end{align*}
\end{proposition}

\begin{proof}
    With a suitable ordering of the vertices of $ G_1 \boxtimes  G_2$ , we get

\begin{align*}
A_\alpha\left( G_1 \boxtimes  G_2\right)=&\left(\begin{array}{ccc}
2\alpha r_1 I_{n_1}+(1-\alpha)A_1 & (1-\alpha)A_1 & O \\
(1-\alpha)A_1 & \alpha (r_1+n_2)I_{n_1} & (1-\alpha)I_{n_1}\otimes J_{1\times n_2} \\
O  & (1-\alpha)I_{n_1}\otimes J_{n_2 \times 1} & I_{n_1}\otimes (\alpha I_{n_2}+A_{\alpha_2})
\end{array}\right).
\end{align*}
The result follows by similar arguments as in Proposition \ref{splitcor}
\end{proof}

Now, in the following corollary, we obtain the $A_\alpha$-eigenvalues of $ G_1 \boxtimes  G_2$, where both $ G_1$ and $ G_2$ are regular graphs.

\begin{corollary}
    Let $ G_i$ be $r_i$-regular graph with $n_i$ vertices for $i=1,2$. Then the $A_\alpha$-spectrum of $ G_1 \boxtimes  G_2$ consists of:
    \begin{enumerate}
        \item $\alpha+\lambda_i(A_{\alpha_2})$ repeated $n_1$ times for each $i=2,3,\dots,n_2$,
        \item two roots of the equation
       $(\lambda-2\alpha r_1-(1-\alpha) \lambda_i(A_1))((\lambda-\alpha-r_2)(\lambda-(r_1+n_2)\alpha)-(1-\alpha)^2n_2)-(1-\alpha)^2(\lambda-\alpha-r_2)\lambda_i(A_1)^2=0,$ for $i=1,2,\dots,n_1$.
    \end{enumerate}
\end{corollary}

Finally, as an application, we construct new pairs of $A_\alpha$-cospectral graphs in the following corollary.

\begin{corollary}
    \begin{enumerate}
    \item Let $ G_1$ and $ G_2$ be two $A$-cospectral regular graphs and $H$ be an arbitrary graph. Then the graphs $ G_1\boxtimes H$ and $ G_2\boxtimes H$ are $A_\alpha$-cospectral.
    \item Let $H_1$ and $H_2$ be two $A_\alpha$-cospectral graphs with $\Gamma_{A_\alpha(H_1)}(\lambda) = \Gamma_{A_\alpha(H_2)}(\lambda)$ for $\alpha \in [0, 1]$. If $ G$ is a regular graph, then the graphs $ G\boxtimes H_1$ and $ G\boxtimes H_2$ are $A_\alpha$-cospectral.
\end{enumerate}
\end{corollary}

Next, we compute $A_\alpha$-spectrum of splitting neighbourhood corona of $ G_1$ and $ G_2$, where $ G_1$ is an $r_1$-reugular graph and $ G_2$ is an arbitrary graph.

\begin{proposition}
     Let $ G_1$ be an $r_1$-regular graph on $n_1$ vertices and $ G_2$ be an arbitrary graph on $n_2$ vertices. Then
     \begin{align*}
         \phi(A_\alpha  \left( G_1 \divideontimes  G_2\right),\lambda) =&\prod_{i=1}^{n_2}\left(\lambda-\alpha r_1-\lambda_i(A_{\alpha_2})\right)^{n_1} \\
&\hspace{-2cm}\prod_{i=1}^{n_1}(\lambda-\alpha r_1)(\lambda-(2+n_2)r_1\alpha-(1-\alpha)^2\Gamma_{A_{\alpha_2}}(\lambda-\alpha r_1)\lambda_i(A_1)^2-(1-\alpha)\lambda_i(A_1))-(1-\alpha)^2\lambda_i(A_1)^2.
     \end{align*}
\end{proposition}

\begin{proof}
    With a suitable ordering of the vertices of $ G_1 \divideontimes  G_2$ , we get

\begin{align*}
A_\alpha\left( G_1 \divideontimes  G_2\right)=&\left(\begin{array}{ccc}
(n_2+2)r_1\alpha I_{n_1}+(1-\alpha)A_1 & (1-\alpha)A_1 & (1-\alpha)J_{1\times n_2}\otimes A_1 \\
(1-\alpha)A_1 & \alpha r_1I_{n_1} & O \\
(1-\alpha)J_{n_2\times1}\otimes A_1  & O & (\alpha r_1 I_{n_2}+A_{\alpha_2})\otimes I_{n_1}
\end{array}\right).
\end{align*}
The result follows by similar arguments as in Proposition \ref{splitcor}
\end{proof}

Now, in the following corollary, we obtain the $A_\alpha$-eigenvalues of $ G_1 \divideontimes  G_2$, where both $ G_1$ and $ G_2$ are regular graphs.

\begin{corollary}
    Let $ G_i$ be $r_i$-regular graph with $n_i$ vertices for $i=1,2$. Then the $A_\alpha$-spectrum of $ G_1 \divideontimes  G_2$ consists of:
    \begin{enumerate}
        \item $\alpha r_1+\lambda_i(A_{\alpha_2})$ repeated $n_1$ times for each $i=2,3,\dots,n_2$,
        \item two roots of the equation
       $(\lambda-\alpha r_1)((\lambda-\alpha r_1-r_2)(\lambda-(2+n_2)r_1\alpha-(1-\alpha)\lambda_i(A_1))-(1-\alpha)^2n_2\lambda_i(A_1)^2)-(1-\alpha)^2(\lambda-\alpha r_1-r_2)\lambda_i(A_1)^2=0,$ for $i=1,2,\dots,n_1$.
    \end{enumerate}
\end{corollary}

Finally, as an application, we construct new pairs of $A_\alpha$-cospectral graphs in the following corollary.

\begin{corollary}
    \begin{enumerate}
    \item Let $ G_1$ and $ G_2$ be two $A$-cospectral regular graphs and $H$ be an arbitrary graph. Then the graphs $ G_1\divideontimes H$ and $ G_2\divideontimes H$ are $A_\alpha$-cospectral.
    \item Let $H_1$ and $H_2$ be two $A_\alpha$-cospectral graphs with $\Gamma_{A_\alpha(H_1)}(\lambda) = \Gamma_{A_\alpha(H_2)}(\lambda)$ for $\alpha \in [0, 1]$. If $ G$ is a regular graph, then the graphs $ G\divideontimes H_1$ and $ G\divideontimes H_2$ are $A_\alpha$-cospectral.
\end{enumerate}
\end{corollary}

\section{\texorpdfstring{$A_\alpha$}{A alpha}-spectrum of corona based on \texorpdfstring{$Q$}{} graph}\label{qcor}
In this section, we compute of $A_\alpha$-spectrum of $ G_1\sqcupdot  G_2$ and $ G_1\sqcupplus  G_2$of the graphs $ G_1$ and $ G_2 $. First, we obtain an expression for the $A_\alpha$-characteristic polynomial of these operations, where $ G_1$ is an $r_1$-reugular graph and $ G_2$ is an arbitrary graph.
\begin{proposition}\label{qvertexcorona}
   Let $ G_1$ be an $r_1$-regular graph on $n_1$ vertices and $ G_2$ be an arbitrary graph on $n_2$ vertices. Then
\begin{align*}
\phi(A_\alpha  \left( G_1 \sqcupdot  G_2\right),\lambda) =&(\lambda-2\alpha r_1+2(1-\alpha))^{m_1-n_1}\prod_{i=1}^{n_2}\left(\lambda-\alpha-\lambda_i(A_{\alpha_2})\right)^{n_1} \\
&\hspace{=-1.5cm}\prod_{i=1}^{n_1}\left((\lambda-\alpha( r_1+n_2)-(1-\alpha)^2\Gamma_{A_{\alpha_2}}(\lambda-\alpha))(\lambda-2\alpha r_1-(1-\alpha)(\lambda_i(A_1)+r_1-2))-(1-\alpha)^2(\lambda_i(A_1)+r_1)\right).
\end{align*} 
\end{proposition}

\begin{proof}
With a suitable ordering of the vertices of $ G_1 \sqcupdot  G_2$ , we get

\begin{align*}
A_\alpha\left( G_1 \sqcupdot  G_2\right)=&\left(\begin{array}{ccc}
\alpha (r_1+n_2) I_{n_1} & (1-\alpha)R_1 & (1-\alpha)I_{n_1}\otimes J_{1 \times n_2}  \\
(1-\alpha)R_1^T & 2\alpha r_1I_{m_1}+(1-\alpha)B_1 & O \\
(1-\alpha)I_{n_1}\otimes J_{n_2 \times 1}  & O & I_{n_1}\otimes(\alpha I_{n_2}+A_{\alpha_2})
\end{array}\right).\\
\left|\lambda I-A_\alpha\left( G_1 \sqcupdot  G_2\right)\right|=&\left|\begin{array}{ccc}
(\lambda-\alpha (r_1+n_2)) I_{n_1} & -(1-\alpha)R_1 & -(1-\alpha)I_{n_1}\otimes J_{1 \times n_2}  \\
-(1-\alpha)R_1^T & (\lambda-2\alpha r_1)I_{m_1}-(1-\alpha)B_1 & O \\
-(1-\alpha)I_{n_1}\otimes J_{n_2 \times 1}  & O & I_{n_1}\otimes((\lambda-\alpha )I_{n_2}-A_{\alpha_2})
\end{array}\right|\\
=&\prod_{i=1}^{n_2}(\lambda-\alpha-\lambda_i(A_{\alpha_2}))^{n_1}\det S\\
\intertext{where }
\det S=&\left|\begin{array}{cc}
    (\lambda-\alpha (r_1+n_2) -(1-\alpha)^2\Gamma_{A_{\alpha_2}}(\lambda-\alpha))I &-(1-\alpha)R_1\\
    -(1-\alpha)R_1^T & (\lambda-2\alpha r_1)I_{m_1}-(1-\alpha)B_1
\end{array}\right|\\
=&(\lambda-\alpha(r_1+n_2)-(1-\alpha)^2\Gamma_{A_{\alpha_2}}(\lambda-\alpha))^{n_1}\\
&\left|(\lambda-2\alpha r_1)I_{m_1}-(1-\alpha)B_1-\frac{(1-\alpha)^2(B_1+2I_{m_1})}{\lambda-\alpha(r_1+n_2)-(1-\alpha)^2\Gamma_{A_{\alpha_2}}(\lambda-\alpha)}\right|\\
\det S=&(\lambda-2\alpha r_1+2(1-\alpha))^{m_1-n_1}\\
&\hspace{=-1.5cm}\prod_{i=1}^{n_1}\left((\lambda-\alpha( r_1+n_2)-(1-\alpha)^2\Gamma_{A_{\alpha_2}}(\lambda-\alpha))(\lambda-2\alpha r_1-(1-\alpha)(\lambda_i(A_1)+r_1-2))-(1-\alpha)^2(\lambda_i(A_1)+r_1)\right).
\intertext{Then the result follows.}
\end{align*}
\end{proof}

Now, in the following corollary, we obtain the $A_\alpha$-eigenvalues of $ G_1 \sqcupdot  G_2$, where both $ G_1$ and $ G_2$ are regular graphs.

\begin{corollary}
    Let $ G_i$ be $r_i$-regular graph with $n_i$ vertices for $i=1,2$. Then the $A_\alpha$-spectrum of $ G_1 \sqcupdot  G_2$ consists of:
    \begin{enumerate}
    \item $2\alpha r_1-2(1-\alpha)$ repeated $m_1-n_1$ times,
        \item $\alpha+\lambda_i(A_{\alpha_2})$ repeated $n_1$ times for each $i=2,3,\dots,n_2$,
        \item three roots of the equation
        $$((\lambda-\alpha-r_2)(\lambda-\alpha(r_1+n_2)-(1-\alpha)^2n_2)(\lambda-2\alpha r_1-(1-\alpha)(\lambda_i(A_1)+r_1-2))-(1-\alpha)^2(\lambda-\alpha-r_2)(\lambda_i(A_1)+r_1)=0,$$ for $i=1,2,\dots,n_1$.
    \end{enumerate}
\end{corollary}

Finally, as an application, we construct new pairs of $A_\alpha$-cospectral graphs in the following corollary.

\begin{corollary}
    \begin{enumerate}
    \item Let $ G_1$ and $ G_2$ be two $A$-cospectral regular graphs and $H$ be an arbitrary graph. Then the graphs $ G_1\sqcupdot H$ and $ G_2\sqcupdot H$ are $A_\alpha$-cospectral.
    \item Let $H_1$ and $H_2$ be two $A_\alpha$-cospectral graphs with $\Gamma_{A_\alpha(H_1)}(\lambda) = \Gamma_{A_\alpha(H_2)}(\lambda)$ for $\alpha \in [0, 1]$. If $ G$ is a regular graph, then the graphs $ G\sqcupdot H_1$ and $ G\sqcupdot H_2$ are $A_\alpha$-cospectral.
\end{enumerate}
\end{corollary}

Next, we compute $A_\alpha$-spectrum of $ G_1\sqcupplus  G_2$, where $ G_1$ is an $r_1$-regular graph and $ G_2$ is an arbitrary graph.
\begin{proposition}
   Let $ G_1$ be an $r_1$-regular graph on $n_1$ vertices and $ G_2$ be an arbitrary graph on $n_2$ vertices. Then
\begin{align*}
\phi(A_\alpha  \left( G_1 \sqcupplus  G_2\right),\lambda) =&(\lambda-\alpha r_1)^{n_1-m_1}\prod_{i=1}^{n_2}\left(\lambda-\alpha-\lambda_i(A_{\alpha_2})\right)^{m_1} \\
&\hspace{=-1cm}((\lambda-\alpha r_1)(\lambda-\alpha(2r_1+n_2)-(1-\alpha)^2\Gamma_{A_{\alpha_2}}(\lambda-\alpha))-2(1-\alpha)^2+2(1-\alpha)(\lambda-\alpha r_1+1-\alpha))^{m_1-n_1}\\
&\hspace{=-3cm}\prod_{i=1}^{n_1}\left((\lambda-\alpha r_1)(\lambda-\alpha(2r_1+n_2)-(1-\alpha)^2\Gamma_{A_{\alpha_2}}(\lambda-\alpha))-2(1-\alpha)^2-(1-\alpha)(\lambda-\alpha r_1+1-\alpha)(\lambda_i(A_1)+r_1-2)\right).
\end{align*} 
\end{proposition}
\begin{proof}
With a suitable ordering of the vertices of $ G_1 \sqcupplus  G_2$ , we get
$$A_\alpha\left( G_1 \sqcupplus  G_2\right)=\left(\begin{array}{ccc}
\alpha r_1 I_{n_1} & (1-\alpha)R_1 & O  \\
(1-\alpha)R_1^T &\alpha(2r_1+n_2)I_{m_1}+(1-\alpha)B_1 & (1-\alpha)I_{m_1}\otimes J_{1 \times n_2} \\
O  & (1-\alpha)I_{m_1}\otimes J_{n_2 \times 1} & I_{m_1}\otimes(\alpha I_{n_2}+A_{\alpha_2})
\end{array}\right).$$
    The proof is similar to the proof of Proposition \ref{qvertexcorona}
\end{proof}

Now, in the following corollary, we obtain the $A_\alpha$-eigenvalues of $ G_1 \sqcupplus  G_2$, where both $ G_1$ and $ G_2$ are regular graphs.

\begin{corollary}
    Let $ G_i$ be $r_i$-regular graph with $n_i$ vertices for $i=1,2$. Then the $A_\alpha$-spectrum of $ G_1 \sqcupplus  G_2$ consists of:
    \begin{enumerate}
    \item $\alpha r_1$ repeated $n_1-m_1$ times,
        \item $\alpha+\lambda_i(A_{\alpha_2})$ repeated $m_1$ times for each $i=2,3,\dots,n_2$,
        \item three roots of the equation
        $$(\lambda-\alpha r_1)((\lambda-\alpha -r_2)(\lambda-\alpha(2r_1 +n_2)-(1-\alpha)^2n_2)-(\lambda-\alpha- r_2)(2(1-\alpha)^2-2(1-\alpha)(\lambda-\alpha r_1+1-\alpha))=0,$$  repeated $m_1-n_1$,
        \item three roots of the equation
        $$(\lambda-\alpha r_1)((\lambda-\alpha -r_2)(\lambda-\alpha(2r_1 +n_2)-(1-\alpha)^2n_2)-(\lambda-\alpha- r_2)(2(1-\alpha)^2+(1-\alpha)(\lambda-\alpha r_1+1-\alpha)(\lambda_i(A_1)+r_1-2))=0,$$  for $i=1,2,\dots,n_1$.
    \end{enumerate}
\end{corollary}

Finally, as an application, we construct new pairs of $A_\alpha$-cospectral graphs in the following corollary.

\begin{corollary}
    \begin{enumerate}
    \item Let $ G_1$ and $ G_2$ be two $A$-cospectral regular graphs and $H$ be an arbitrary graph. Then the graphs $ G_1\sqcupplus H$ and $ G_2\sqcupplus H$ are $A_\alpha$-cospectral.
    \item Let $H_1$ and $H_2$ be two $A_\alpha$-cospectral graphs with $\Gamma_{A_\alpha(H_1)}(\lambda) = \Gamma_{A_\alpha(H_2)}(\lambda)$ for $\alpha \in [0, 1]$. If $ G$ is a regular graph, then the graphs $ G\sqcupplus H_1$ and $ G\sqcupplus H_2$ are $A_\alpha$-cospectral.
\end{enumerate}
\end{corollary}

\bibliographystyle{unsrt}  
\bibliography{references}

\end{document}